\documentclass{article}
\usepackage {amsmath, amsfonts, amssymb, mathrsfs, stmaryrd, amsthm, diagrams,accents, sectsty, txfonts, ifthen,calc}
\usepackage[usenames]{color}
\def\mat#1#2#3#4{\left(\begin{array}{cc} #1 & #2 \\ #3 & #4 \end{array}\right)}
\def\smat#1#2#3#4{{\tiny\mat{#1}{#2}{#3}{#4}}}
\def\vthree#1#2#3{\left(\begin{array}{c}#1\\#2\\#3\end{array}\right)}

\def\mb#1{\mathbb{#1}}
\def\tx#1{{\rm #1}}
\def\tb#1{\textbf{#1}}
\def\ti#1{\textit{#1}}

\def\R{\mathbb{R}}
\def\C{\mathbb{C}}
\def\Q{\mathbb{Q}}

\def\Z{\mathbb{Z}}

\def\ul#1{\underline{#1}}

\def\hat{\widehat}

\newarrow{Equals}{=}{=}{=}{=}{}

\newtheorem{thm}{Theorem}[section]
\newtheorem{lem}[thm]{Lemma}
\newtheorem{pro}[thm]{Proposition}
\newtheorem{cor}[thm]{Corollary}
\newtheorem{dfn}[thm]{Definition}
\newtheorem{fct}[thm]{Fact}

\def\rmk{\tb{Remark: }}

\def\pf{\tb{Proof: }}
\def\nt{\tb{Note: }}

\newenvironment{mytitle}
{\begin{center}\large\sc}
{\end{center}}

\begin{document}

\begin{mytitle} Decomposition of splitting invariants in split real groups \end{mytitle}
\begin{center} Tasho Kaletha \end{center}

{\renewcommand{\thefootnote}{}}

\begin{abstract}
To a maximal torus in a quasi-split semi-simple simply-connected group over a local field of characteristic 0, Langlands and Shelstad construct in \cite{LS1} a
cohomological invariant called the splitting invariant, which is an important component of their endoscopic transfer factors. We study this invariant in the
case of a split real group and prove a decomposition theorem which expresses this invariant for a general torus as a product of the corresponding invariants
for simple tori. We also show how this reduction formula allows for the comparison of splitting invariants between different tori in the given real group.
\end{abstract}

In applications of harmonic analysis and representation theory of reductive groups over local fields to questions in number theory, a central role is played by
the theory of endoscopy. This theory associates to a given connected reductive group $G$ over a local field $F$ a collection of connected reductive groups over
$F$, often denoted by $H$, which have smaller dimension (except when $H=G$), but are usually not subgroups of $G$. The geometric side of the theory is then
concerned with transferring functions on $G(F)$ to functions on $H(F)$ in such a way that suitable linear combinations of their orbital integrals are
comparable, while the spectral side is concerned with transferring ``packets'' of representations on $H(F)$ to ``packets'' of representations on $G(F)$ in such
a way that suitable linear combinations of their characters are comparable. In both cases, the comparison involves certain normalizing factors, called
geometric or spectral transfer factors.

Over the real numbers, the theory of endoscopy was developed in a series of profound papers by Diana Shelstad, in which she defines geometric and spectral
transfer factors and proves that these factors indeed give a comparison of orbital integrals and character formulas between $G$ and $H$. A very subtle and
complicated feature of the transfer factors was the need to assign a $\pm$-sign to each maximal torus in $G$ in a coherent manner, and Shelstad was able to
prove that this is possible. A uniform and explicit definition of geometric transfer factors for all local fields was given in \cite{LS1}. An explicit
construction of spectral transfer factors over the real numbers was given in \cite{S2}, while over the $p$-adic numbers their existence is still conjectural.
The structure of transfer factors is quite complex -- both the geometric and the real spectral ones are a product of multiple terms of group-theoretic or
Galois-cohomological nature. There are numerous choices involved in the construction of each individual term, but the product is independent of most choices.
One term that is common for both the geometric and the real spectral transfer factors is called $\Delta_I$. It is regarded as the most subtle and is the one
that makes explicit the choice of coherent collection of signs in Shelstad's earlier work. At its heart is a Galois-cohomological object, called the splitting
invariant. The splitting invariant is an element of $H^1(F,T)$ associated to any maximal torus $T$ of a quasi-split semi-simple simply-connected group $G$,
whose construction occupies the first half of \cite[Sec.2]{LS1}. It is depends on the choice of a splitting $(T_0,B_0,\{X_\alpha\}_{\alpha \in \Delta})$ of $G$
as well as a-data $\{a_\beta\}_{\beta\in R(T,G)}$.

This paper addresses the following question: If one has two tori in a given real group which originate from the same endoscopic group, how can one compare
their splitting invariants? While there will in general be no direct relation between $H^1(F,T_1)$ and $H^1(F,T_2)$ for two tori $T_1$ and $T_2$ of $G$, if
both those tori originate from $H$ then there are certain natural quotients of their cohomology groups which are comparable, and it is the image of the
splitting invariant in those quotients that is relevant to the construction of $\Delta_I$. An example of a situation where this problem arises is the
stabilization of the topological trace formula of Goresky-MacPherson. One is led to consider characters of virtual representations which occur as sums indexed
over tori in $G$ that originate from the same endoscopic group $H$, and each summand carries a $\Delta_I$-factor associated to the corresponding torus.

To describe the results of this paper, let $G$ be a split simply-connected real group and $(T_0,B_0,\{X_\alpha\}_{\alpha \in \Delta})$ be a fixed splitting.
For a subset $A \subset R$, consisting of strongly orthogonal roots, let $S_A$ denote the element of the Weyl group of $T_0$ given by the product of the
reflections associated to the elements of $A$ (the order in which the product is taken is irrelevant). We show that associated to $A$ there is a canonical
maximal torus $T_A$ of $G$ and a set of isomorphisms of real tori $T_0^{S_A} \rightarrow T_A$, where $T_0^{S_A}$ is the twist by $S_A$ of $T_0$. Any maximal
torus in $G$ is $G(\R)$-conjugate to one of the $T_A$, so it is enough to study the tori $T_A$. We give an expression in purely root-theoretic terms for a
certain 1-cocycle in $Z^1(\R,T_0^{S_A})$. This cocycle has the property that its image in $Z^1(\R,T_A)$ under any of the isomorphisms $T^{S_A} \rightarrow T_A$
above is the same, and the class in $H^1(\R,T_A)$ of that image is the splitting invariant of $T_A$ (associated to a specific choice of a-data). Moreover, we
prove a reduction theorem which shows that this cocycle is a product over $\alpha \in A$ of the cocycles associated to the canonical tori $T_{\{\alpha\}}$,
thereby reducing the study of the splitting invariant of $T_A$ to those of the various $T_{\{\alpha\}}$. This product decomposition takes place inside the
group $Z^1(\R,T_0^{S_A})$ -- that is, we show that the elements of $Z^1(\R,T_0^{s_\alpha})$ associated to the various $T_{\{\alpha\}}$ with $\alpha \in A$ also
lie in $Z^1(\R,T_0^{S_A})$ and that their product is the element associated to $T_A$. Finally we show that if $A' \subset A$ and the tori $T_{A'}$ and $T_A$
originate from the same endoscopic group, then the endoscopic characters on the cohomology groups $H^1(\R,T_{A'})$ and $H^1(\R,T_A)$ factor through certain
explicitly given quotients of these groups, and the quotient of $H^1(\R,T_{A'})$ is canonically embedded into that of $H^1(\R,T_A)$. This, together with the
reduction theorem, allows for a direct comparison of the values that the endoscopic characters associate to the splitting invariants for $T_{A'}$ and $T_A$.

The paper is organized as follows: Section \ref{sec:prelims} contains a few basic facts and serves mainly to fix notation for the rest of the paper. Section
\ref{sec:strong} contains proofs of general facts about subsets of strongly orthogonal roots in reduced root systems, which are needed as a preparation for the
reduction theorem mentioned above. The study of the splitting invariants takes place in section \ref{sec:splitting}, where first the splitting invariant for
the tori $T_{\{\alpha\}}$ is computed, and after that the results of section \ref{sec:strong} are used to reduce the case of $T_A$ to that of $T_{\{\alpha\}}$.
While the statement of the reduction theorem appears natural and clear, the proof contains some subtle points. First, one has to choose the Borel $B_0$ in the
splitting of $G$ with care according to the strongly orthogonal set $A$. As remarked in section \ref{sec:splitting}, this choice does not affect the splitting
invariant, but it significantly affects its computation. Moreover, the root system $G_2$ exhibits a singular behavior among all reduced root systems as far as
pairs of strongly-orthogonal roots are concerned. Section \ref{sec:explicit} contains explicit computations of the splitting invariants of the tori
$T_{\{\alpha\}}$ for all split almost-simple classical groups. In section \ref{sec:comparison} we construct the aforementioned quotients of the cohomology
groups and the embedding between them.

%We'd like to thank Professor Robert Kottwitz and Professor Diana Shelstad for their valuable comments.

%---------------------------------------------------------------------------------
% Section
%=================================================================================

\section{Notation and preliminaries}\label{sec:prelims}

Throughout this paper $G$ will stand for a split semi-simple simply-connected group over $\mathbb{R}$ and $(B_0, T_0, \{X_\alpha\})$ will be a splitting of
$G$. We write $R=R(T_0,G)$ for the set of roots of $T_0$ in $G$, set $\alpha>0$ if $\alpha\in R(T_0,B_0)$, denote by $\Delta$ the set of simple roots in
$R(T_0,B_0)$ and by $\Omega$ the Weyl-group of $R$, which is identified with $N(T_0)/T_0$. Moreover we put $\Gamma=\tx{Gal}(\C/\R)$ and denote by $\sigma$ both
the non-trivial element in that group, as well as its action on $T_0$. The notation $g \in G$ will be shorthand for $g \in G(\C)$, and $\tx{Int}(g)h =
ghg^{-1}$.

\subsection{$sl_2$-triples}

For any $\alpha \in R(T_0,B_0)$ we have the coroot $\alpha^\vee: \mathbb{G}_m \rightarrow T_0$ and its differential $d\alpha^{\vee}: \mathbb{G}_a \rightarrow
\tx{Lie}(T_0)$. We put $H_\alpha := d\alpha^{\vee}(1) \in \tx{Lie}(T_0)$. Given $X_\alpha \in \tx{Lie}(G)_\alpha$ non-zero, there exists a unique $X_{-\alpha}
\in \tx{Lie}(G)_{-\alpha}$ so that $[H_\alpha, X_\alpha, X_{-\alpha}]$ is an $sl_2$-triple. The map
\[ \mat{1}{0}{0}{-1} \mapsto H_\alpha,\quad \mat{0}{1}{0}{0} \mapsto X_\alpha,\quad \mat{0}{0}{1}{0} \mapsto X_{-\alpha} \]

gives a homomorphism $sl_2 \rightarrow \tx{Lie}(G)$ which integrates to a homomorphism $\tx{SL}_2 \rightarrow G$ and one has

\begin{diagram}
sl_2 & \rTo & \tx{Lie}(G) \\
\dTo^{\tx{exp}} & & \dTo_{\tx{exp}}\\
\tx{SL}_2 & \rTo & G&\\
\end{diagram}

The image of $\smat{a}{b}{c}{d} \in \tx{SL}_2$ under this homomorphism will be called $\smat{a}{b}{d}{c}_{X_\alpha}$. Notice that
$\smat{t}{0}{0}{t^{-1}}_{X_\alpha} = \alpha^\vee(t)$.

\begin{fct} \label{fct:strongcommute}
Let $\alpha,\beta \in R$ be s.t. $\alpha+\beta \notin R$ and $\alpha-\beta \notin R$. For any non-zero elements $X_\alpha \in \tx{Lie}(G)_\alpha$ and $X_\beta
\in \tx{Lie}(G)_\beta$, the homomorphisms $\varphi_{X_\alpha}, \varphi_{X_\beta} : \tx{SL}_2 \rightarrow G$ given by $X_\alpha$ and $X_\beta$ commute.
\end{fct} \pf Since for any field $k$, $\tx{SL}_2(k)$ is generated by its two subgroups
\[ \left\{ \mat{1}{u}{0}{1}\left|\right.\ u \in k\right\} \qquad \left\{ \mat{1}{0}{u}{1}\left|\right.\ u \in k\right\} \]
it is enough to show that, for any $u,v \in \C$, each of $\tx{exp}(uX_\alpha)$ and $\tx{exp}(uX_{-\alpha})$ commutes with each of $\tx{exp}(vX_\beta)$ and
$\tx{exp}(vX_{-\beta})$. This follows from \cite[10.1.4]{Spr} and our assumption on $\alpha,\beta$.\qed

\subsection{Chevalley bases}\label{subsec:chevbas}

For $\alpha \in \Delta$ let $n_\alpha = \tx{exp}(X_\alpha) \tx{exp}(-X_{-\alpha}) \tx{exp}(X_\alpha) = \smat{0}{1}{-1}{0}_{X_\alpha}$. Given $\mu \in \Omega$
we have the lift $n(\mu) \in N(T_0)$ given by
\[n(\mu) = n_{\alpha_1} \cdots n_{\alpha_q} \]
where $s_{\alpha_1} \cdots s_{\alpha_q} = \mu$ is any reduced expression (by \cite[11.2.9]{Spr} this lift is independent of the choice of reduced expression).
Notice $n(\mu) \in N(T_0)(\mb{R})$ since $T_0$ is split. Put
\[ X_{\mu|\alpha} := \tx{Int}(n(\mu))\cdot X_\alpha \]

Then $X_{\mu|\alpha} \in \tx{Lie}(G)_{\mu\alpha}$ is a non-zero element.

\begin{lem}
If $\alpha, \alpha' \in \Delta$ and $\mu, \mu' \in \Omega$ are s.t. $\mu\alpha = \mu'\alpha'$ then we have in $\tx{Lie}(G)_{\mu\alpha}$ the equality
\[X_{\mu'|\alpha'} = \prod_{\substack{\beta>0\\(\mu')^{-1} \beta<0\\\mu^{-1} \beta > 0}} (-1)^{\langle \beta^\vee, \mu\alpha \rangle} \cdot X_{\mu|\alpha} \]
\end{lem}

\pf By \cite[11.2.11]{Spr} the relation $(\mu')^{-1} \cdot \mu \alpha = \alpha'$ implies
\[X_{\alpha'} = \tx{Int}\left[ n\left( \left( \mu'\right)^{-1} \cdot \mu\ \right) \right] X_\alpha \]

The claim now follows from \cite[2.1.A]{LS1} and the following computation
\begin{eqnarray*}
X_{\mu'|\alpha'} & = & \tx{Int}\left( n \left(\mu' \right)\right)X_{\alpha'}\\
& = & \tx{Int} \left[ n\left(\mu'\right)n\left(\left(\mu'\right)^{-1} \mu\right) \right] X_\alpha\\
& = & \tx{Int} \left[ t\left(\mu', (\mu')^{-1}\mu \right) \cdot n(\mu) \right] X_\alpha\\
& = & \tx{Int} \left[ t\left(\mu', (\mu')^{-1} \mu \right) \right]X_{\mu|\alpha}\\
& = & (\mu \alpha) \left(t(\mu', (\mu')^{-1}\mu)\right) \cdot X_{\mu|\alpha}
\end{eqnarray*}

\qed

\rmk We see that while the "absolute value" of $X_{\mu|\alpha}$ only depends on the root $\mu \cdot \alpha$, its "sign" does depend on both $\mu$ and $\alpha$.

\begin{dfn} For $\gamma \in R$, $\mu, \mu' \in \Omega$ put
\[ \epsilon (\mu', \gamma, \mu) := \prod_{\substack{\beta>0\\(\mu')^{-1} \beta<0\\\mu^{-1} \beta > 0}} (-1)^{\langle \beta^\vee, \gamma \rangle}\]
\end{dfn}
\rmk With this definition we can reformulate the above lemma as follows
\begin{cor}\label{cor:chevtrans}
If $\gamma \in R$ and $\mu, \mu' \in \Omega$ are s.t. $\mu^{-1}\gamma, (\mu')^{-1}\gamma \in \Delta$ then
\[ X_{\mu'|(\mu')^{-1}\gamma} = \epsilon(\mu', \gamma, \mu) \cdot X_{\mu|\mu^{-1}\gamma} \]
\end{cor}

\rmk If  for each $\gamma \in R$ we choose $\mu_\gamma \in \Omega$ so that $\mu_\gamma^{-1} \gamma \in \Delta$, then $\{ X_{\mu_\gamma |
\mu_\gamma^{-1}\gamma}\}_{\gamma \in R}$ is a Chevalley system in the sense of \cite[exp XXIII \S6]{SGAIII.3}.

\subsection{Cayley-transforms}\label{sec:cayley}

Let $\alpha \in R(T_0, B_0)$ and choose $X_\alpha \in \tx{Lie}(G)_\alpha(\mathbb{R})-\{0\}$. Put
\[ g_\alpha := \tx{exp}\left( \frac{i\pi}{4} (X_\alpha + X_{-\alpha}) \right)\]
Then $\sigma(g_\alpha) = \tx{exp}\left( -\frac{i\pi}{4} (X_\alpha + X_{-\alpha}) \right) = g_\alpha^{-1}$ and $\sigma(g_\alpha)^{-1} \cdot g_\alpha =
g_\alpha^2 = \tx{exp}\left( \frac{i\pi}{2} (X_\alpha + X_{-\alpha}) \right)$. We have
\[ g_\alpha = \left[ \frac{\sqrt{2}}{2} \mat{1}{i}{i}{1} \right]_{X_\alpha}, \quad\
   g_\alpha^2 = \mat{0}{i}{i}{0}_{X_\alpha}, \quad g_\alpha^4 = \mat{-1}{0}{0}{-1}_{X_\alpha} = \alpha^\vee(-1) \]

\begin{fct} The images of $T_0$ under $\tx{Int}(g_\alpha)$ and $\tx{Int}(g_\alpha^{-1})$ are the same. They are a torus $T$ defined over $\mathbb{R}$
and the transports of the $\Gamma$-action on $T$ to $T_0$ via $\tx{Int}(g_\alpha^{-1})$ and $\tx{Int}(g_\alpha)$ both equal $s_\alpha \rtimes \sigma$.
\end{fct}

\pf
\[ \tx{Int}(g_\alpha) T_0 = \tx{Int}(g_\alpha^{-1}) \tx{Int}(g_\alpha^2) T_0 = \tx{Int}(g_\alpha^{-1}) s_\alpha T_0 = \tx{Int}(g_\alpha^{-1}) T_0 \]
\[ \sigma(\tx{Int}(g_\alpha)T_0) = \tx{Int}(\sigma(g_\alpha))T_0 = \tx{Int}(g_\alpha^{-1}) T_0 \]
\[ \tx{Int}(\sigma(g_\alpha)^{-1} g_\alpha) = \tx{Int}(g_\alpha^2) = s_\alpha = \tx{Int}(g_\alpha^{-2}) = \tx{Int}(\sigma(g_\alpha) g_\alpha^{-1}) \]
\qed

\nt Different choices of $X_\alpha$ will lead to different (yet conjugate) tori $T$. However, since we have fixed a splitting there is up to a sign a canonical
$X_\alpha$. Changing the sign of $X_\alpha$ changes $g_\alpha$ to $g_\alpha^{-1}$, hence $T$ does not change. Thus we conclude:

The choice of a splitting gives for each $\alpha \in R(T_0, B_0)$ the following canonical data:
\begin{enumerate}
\item a pair $\{X,X'\} \subseteq \tx{Lie}(G)_\alpha(\mb{R})-\{0\}$ with $X' = -X$
\item a torus $T_\alpha$ on which $\Gamma$ acts via $s_\alpha \rtimes \sigma$,
\item a pair $\varphi, \varphi'$ of isomorphisms $T_0^{s_\alpha} \rightarrow T_\alpha$ s.t. $\varphi' = \varphi \circ s_\alpha$, given by the Cayley-transforms
with respect to $X,X'$.
\end{enumerate}

\begin{cor} For $\alpha \in R(T_0,B_0)$ let $T_\alpha$ be the canonically given torus as above. For $\mu,\mu' \in \Omega$ s.t.
$\mu^{-1}\alpha,(\mu')^{-1}\alpha \in \Delta$ let $\varphi,\varphi' : T_0^{s_\alpha} \rightarrow T_\alpha$ be the isomorphisms given by
$\tx{Int}(g_{X_{\mu|\mu^{-1}\alpha}})$ and $\tx{Int}(g_{X_{\mu'|(\mu')^{-1}\alpha}})$. Then
\[ \varphi' = \begin{cases} \varphi&,\epsilon(\mu',\alpha,\mu)=1\\ \varphi\circ s_\alpha&, \epsilon(\mu',\alpha,\mu)=-1\end{cases} \]
\end{cor}
\pf Clear.\qed

\tb{Notation:} From now on we will write $g_{\mu,\alpha}$ instead of $g_{X_{\mu|\mu^{-1}\alpha}}$. This notation will only be employed in the case that
$\alpha \in R(T_0,B_0)$ and $\mu^{-1}\alpha \in \Delta$.

%\subsection{Transport of roots}
%
%We work over $\mb{C}$. Let $\varphi : T_0 \rightarrow T$ be an
%isomorphism of maximal tori. The diagrams
%%\begin{diagram}
%%T_0 & \rTo^\varphi & T & \rTo & \mb{G}_m \quad,\\
%%\mb{G}_m & \rTo & T_0 & \rTo^\varphi & T\\
%%\end{diagram}
%
%\[ T_0
%\stackrel{\varphi}{\longrightarrow}T\stackrel{}{\longrightarrow}
%\mb{G}_m\ ,\qquad \mb{G}_m \stackrel{}{\longrightarrow} T_0
%\stackrel{\varphi}{\longrightarrow} T \]
%
%give maps $\varphi^* : R(T,G) \rightarrow R(T_0,G)$, $\varphi_* : R^\vee(T_0, G) \rightarrow
%R^\vee(T,G)$ s.t. for $\alpha \in R(T,G)$, $\lambda \in R^\vee(T_0, G)$ \[ \langle \lambda ,
%\varphi^*\alpha \rangle = \langle \varphi_* \lambda , \alpha \rangle \]
%
%If $\varphi = \tx{Int}(h)$ for $h \in G$ and we put $B = \tx{Int}(h) B_0$, then
%\begin{diagram}
%\mb{G}_a &\rTo^{u_\beta} & B_0 &\rTo^{\tx{Int}(h)}&B\\
%\end{diagram}
%gives \[ \tx{Int}(h)u_\beta = u_{\tx{Int}(h^{-1})^*\beta} \]
%
%If $h \in N(T_0)$ write $\tx{Int}(h)|_{T_0} = \omega \in \Omega(T_0)$ and
%\begin{eqnarray*}
%\omega \cdot \beta = (\omega^{-1})^*\beta = \beta \circ \omega^{-1} \quad,\quad \beta \in R(T_0,
%G)\\
%\omega \cdot \lambda = \omega_*\lambda = \omega \circ \lambda \quad,\quad \lambda \in R^\vee(T_0,
%G)
%\end{eqnarray*}
%Then
%\begin{eqnarray*}
%&&\langle \lambda, \omega\beta \rangle = \langle \omega^{-1}\lambda, \beta \rangle\\
%&& \tx{Int}(h)u_{\beta} = u_{\omega\beta}
%\end{eqnarray*}

%---------------------------------------------------------------------------------
% Section
%=================================================================================

\section{Strongly orthogonal subsets of root systems}\label{sec:strong}

In this section, a few technical facts about strongly orthogonal subsets of root systems are proved.

\begin{dfn}\ \\[-25pt]
\begin{enumerate}
\item $\alpha,\beta \in R$ are called strongly orthogonal if $\alpha+\beta \notin R$ and $\alpha-\beta \notin R$.
\item $A \subset R$ is called a strongly orthogonal subset (SOS) if it consists of pairwise strongly orthogonal roots.
\item $A \subset R$ is called a maximal strongly orthogonal subset (MSOS) if it is a SOS and is not properly contained in a SOS.
\end{enumerate}
\end{dfn}

A classification of the Weyl group orbits of MSOS in irreducible root systems was given in \cite{AK}. In some cases, there exists more than one orbit. To handle
these cases, we will use the following definition and lemma.

\begin{dfn} Let $A_1$,$A_2$ be SOS in $R$. $A_2$ will be called adapted to $A_1$ if $\tx{span}(A_2) \subset \tx{span}(A_1)$ and for all distinct $\alpha,\beta \in A_2$
\[ \{ a \in A_1:\ (a,\alpha) \neq 0 \} \cap \{ a \in A_1:\ (a,\beta) \neq 0 \} = \emptyset \]
where $()$ is any $\Omega$-invariant scalar product on the real vector space spanned by $R$.
\end{dfn}
Note that any $A$ is adapted to itself.

\begin{lem} \label{lem:msosreps} There exist representatives $A_1,...,A_k$ of the Weyl group orbits of MSOS s.t. $A_1$ has maximal length and $A_2,...,A_k$ are adapted to $A_1$. \end{lem}
\pf This follows from the explicit classification in \cite{AK}.\qed

\tb{Notation:} If $A$ is a SOS then all reflections with respect to elements in $A$ commute. Their product will be denoted by $S_A$.

\begin{dfn}
For a root system $R$, a choice of positive roots $>$ and a subset $A$ of $R$ let $\#(R,>,A)$ be the following statement
\begin{eqnarray*}
\forall \alpha_1,\alpha_2 \in A\ \forall \beta>0&&\\
&&\alpha_1\neq\alpha_2\ \wedge\ s_{\alpha_1}(\beta) < 0 \Longrightarrow s_{\alpha_2}(\beta) >0
\end{eqnarray*}
and let $\#\#(R,>,A)$ be the following statement
\begin{eqnarray*}
\forall A_1,A_2 \subset A&& \forall \beta>0 \\
&&A_1 \cap A_2 = \emptyset\ \wedge\ S_{A_1}(\beta) <0 \Longrightarrow S_{A_2}(\beta)>0\ \wedge\ S_{A_1}S_{A_2}(\beta)<0
\end{eqnarray*}
\end{dfn}
\rmk We will soon show that these statements are equivalent. Moreover we will show that for any SOS $A \subset R$ we can choose $>$ so that the triple
$(R,>,A)$ verifies these statements. For this it is more convenient to work with $\#$. For the applications however, we need $\#\#$.

\begin{lem}\label{lem:posroots}
Let $R$ be a reduced root system and $A \subset R$ a SOS. There exists a choice of positive roots $>$ s.t. $\#(R,>,A')$ holds for any $A'$ adapted to $A$.
\end{lem}
\pf Let $V$ denote the real vector space spanned by $R$, and $(\ ,\ )$ be an $\Omega$-invariant scalar product on $V$. The elements of $A$ are orthogonal wrt
$(\ ,\ )$. Extend $A$ to an orthogonal basis $(a_1,...,a_n)$ of $V$. Define the following notion of positivity on $R$
\[ \alpha > 0 \Longleftrightarrow (\alpha,a_{i_0}) > 0\qquad \tx{for}\qquad i_0 = \min\{i:\ (\alpha,a_i)\neq 0\} \]

It is clear from the construction that with this notion $\#(R,>,A')$ is satisfied for any $A'$ adapted to $A$. We just need to check that $>0$ defines a choice
of positive roots, which we will now do.

It is clear that for each $\alpha \in R$ precisely one of $\alpha>0$ or $-\alpha>0$ is true. We will construct $p \in V$ s.t. for all $\alpha \in R$
\[ \alpha>0 \Longleftrightarrow (\alpha,p) > 0 \]
Let
\begin{eqnarray*}
m&=&\min\{|(\alpha,a_i)|:\ \alpha \in R,1\leq i \leq n,(\alpha,a_i)\neq 0\}\\
M&=&\max\{|(\alpha,a_i)|:\ \alpha \in R,1\leq i \leq n\}
\end{eqnarray*}
Construct recursively real numbers $p_1,...,p_n$ s.t.
\[ p_n=1\qquad p_i > \frac{M}{m}\sum_{k>i}p_k \]
and put $p = \sum p_ia_i$. If $\alpha \in R$ is s.t. $\alpha>0$ and $i_0$ is the smallest $i$ s.t. $(\alpha,a_i)\neq 0$ then
\[ (\alpha,p) = \sum_{i=i_0}^n p_i(\alpha,a_i) > mp_{i_0} - M\sum_{k>i_0}p_k > 0 \]
Thus
\[ \alpha > 0 \Longrightarrow (\alpha,p)>0 \]
The converse implication follows formally:
\[ \neg (\alpha >0) \Leftrightarrow -\alpha > 0 \Rightarrow (-\alpha,p) > 0 \Rightarrow \neg( (\alpha,p)>0) \]
\qed

\rmk The truth value of the statement $\#(R,>,A)$ and the notion of being adapted to $A$ are unchanged if one replaces elements of $A$ by their negatives. Thus
we can always assume that the elements of $A$ are positive.

\rmk It is necessary to choose the set of positive roots based on $A$ in order for $\#(R,>,A)$ to be true. An example that $\#(R,>,A)$ may be false is provided
by $V=\R^3,R=D_3$ with positive roots
\[ \vthree{1}{-1}{0},\vthree{1}{0}{-1},\vthree{0}{1}{-1},\vthree{1}{1}{0},\vthree{1}{0}{1},\vthree{0}{1}{1} \]
and
\[ A=\left\{\vthree{1}{0}{-1},\vthree{1}{0}{1}\right\}, \beta=\vthree{1}{-1}{0} \]

\begin{fct} \label{fct:g2sos}
Let $R=G_2$ and $>$ any choice of positive roots. All MSOS $A$ of $R$ lie in the same Weyl-orbit and moreover automatically satisfy $\#(R,>,A)$. Some of these
$A$ contain simple roots.
\end{fct}
\pf This is an immediate observation.\qed

\begin{pro} \label{pro:strong} Let $A \subset R$ be a SOS and $>$ be a choice of positive roots. Then the statements $\#(R,>,A)$ and $\#\#(R,>,A)$ are equivalent.
\end{pro}
\pf First, we show that $\#$ implies the following statement, to be called $\#_1$:
\begin{eqnarray*}
\forall \alpha_1,\alpha_2 \in A\ \forall \beta>0&&\\
&&\alpha_1\neq\alpha_2\ \wedge\ s_{\alpha_1}(\beta) < 0 \Longrightarrow s_{\alpha_2}(\beta)>0\ \wedge\ s_{\alpha_1}s_{\alpha_2}(\beta)<0
\end{eqnarray*}

Let $\alpha_1,\alpha_2 \in A$ and $\beta>0$ be s.t. $s_{\alpha_1}(\beta)<0$. Put $\beta' = -s_{\alpha_1}(\beta)$. Then $\beta'>0$ and
$s_{\alpha_1}(\beta')=-\beta<0$. Then $\#$ implies that $s_{\alpha_2}s_{\alpha_1}(\beta)=-s_{\alpha_2}(\beta')<0$.

Next we show that $\#_1$ implies the following statement, to be called $\#_2$:
\begin{eqnarray*}
\forall \alpha_1\in A\forall A_2\subset A\ \forall \beta>0&&\\
&&\alpha_1\notin A_2\ \wedge\ s_{\alpha_1}(\beta) < 0 \Longrightarrow S_{A_2}(\beta)>0\ \wedge\ s_{\alpha_1}S_{A_2}(\beta)<0
\end{eqnarray*}
We do this by induction of the cardinality of $A_2$, the case of $A_2$ singleton being precisely $\#_1$. Now let $\alpha_1 \in A$, $A_2\subset
A\setminus\{\alpha_1\}$, and $\beta>0$ be s.t. $s_{\alpha_1}(\beta)<0$. Choose $\alpha_2 \in A_2$ and put $\beta':=s_{\alpha_2}(\beta)$. Then by $\#_1$ we have
$\beta'>0$ and $s_{\alpha_1}(\beta')<0$. Applying the inductive hypothesis we obtain $S_{A_2}(\beta) = S_{A_2\setminus\{\alpha_2\}}(\beta')>0$ and
$s_{\alpha_1}S_{A_2}(\beta) = s_{\alpha_1}S_{A_2\setminus\{\alpha_2\}}(\beta')<0$.

Now we show that $\#_2$ implies the following statement, to be called $\flat$:

If $A_2 \subset A$ and $\beta>0$ are s.t. $S_{A_2}(\beta)<0$ then there exists $\alpha_2 \in A_2$ s.t. $s_{\alpha_2}(\beta)<0$.

To see this, let $A_3 \subset A_2$ be a subset of minimal size s.t. $S_{A_3}(\beta)<0$. Take $\alpha_3 \in A_3$ and put $\beta' =
S_{A_3\setminus\{\alpha_3\}}(\beta)$. By minimality of $A_3$ we have $\beta'>0$, and moreover $s_{\alpha_3}(\beta')=S_{A_3}(\beta)<0$. Then $\#_2$ implies that
$s_{\alpha_3}(\beta) = s_{\alpha_3}S_{A_3\setminus\{\alpha_3\}}(\beta')<0$.

Finally we show that $\#_2$ implies the statement $\#\#$. Take $A_1,A_2 \subset A$ s.t. $A_1 \cap A_2 = \emptyset$ and $\beta>0$ s.t. $S_{A_1}(\beta)<0$. By
$\flat$ there exists $\alpha_1 \in A_1$ s.t. $s_{\alpha_1}(\beta)<0$. Since $\alpha_1 \notin A_2$ we get from $\#_2$ that $S_{A_2}(\beta)>0$ and
$S_{A_1}S_{A_2}(\beta)=s_{\alpha_1}S_{A_2}S_{A_1\setminus\{\alpha_1\}}(\beta)<0$.

This shows that $\#$ imples $\#\#$. The converse implication is trivial.\qed

\begin{pro} \label{pro:decomp} For an SOS $A \subset R$, and a choice $>$ of positive roots, let
\[ R^+_A = \{ \beta \in R :\ \beta>0\ \wedge\ S_A\beta < 0 \} \]
Assume that $>$ is chosen so that $\#\#(R,>,A)$ is true. Then if $A',A'' \subset A$ are disjoint, so are $R^+_{A'}$ and $R^+_{A''}$, and $R^+_{A'\cup A''} =
R^+_{A'} \cup R^+_{A''}$. Moreover, the action of $S_{A'}$ on $R$ preserves $R^+_{A''}$.
\end{pro}
\pf This follows immediately. \qed

\begin{cor} \label{cor:decomp} If $A$ is a SOS and $>$ is chosen so that $\#(R,>,A)$ is true then
\[ R^+_A = \coprod_{\alpha \in A} R^+_{\alpha} \]
\end{cor}
\pf Clear.

\begin{lem} Let $R$ be a root system, $V$ the real vector space spanned by it, $Q \subset V$ the root lattice, and $(\ ,\ )$ a Weyl-invariant scalar product on $V$.
If $v \in Q$ is s.t.
\[ |v| \leq \min\{|\alpha|:\ \alpha \in R\} \]
where $|\ |$ is the Euclidian norm arising from $(\ ,)$ then $v \in R$ and the above inequality is an equality. \end{lem} \pf Choose a presentation
\[ v=\sum_{\alpha \in R} n_\alpha \alpha,\qquad n_\alpha \in \Z_{\geq 0} \]
s.t. $\sum_{\alpha} n_\alpha$ is minimal. First we claim that if $\alpha,\beta \in R$ contribute to this sum, then $(\alpha,\beta) \geq 0$. If that were not
the case, then by \cite[Ch.VI,\S1,no.3,Thm.1]{Bou} we have that  $\gamma:=\alpha+\beta \in R \cup \{0\}$ and we can replace the contribution $\alpha+\beta$ in
the sum by $\gamma$, contradicting its minimality. Now, if $\gamma \in R$ is any root contributing to the sum, we get
\[ |v|^2 = \sum_{\alpha,\beta \in R}n_\alpha n_\beta (\alpha,\beta) \geq n_\gamma^2(\gamma,\gamma) \geq (\gamma,\gamma) = |\gamma|^2 \]
with equality precisely when $v=\gamma$. \qed

\begin{lem} \label{lem:strongco} Let $R$ be a root system and $\alpha,\beta \in R$ two strongly orthogonal roots.
If $\alpha^\vee+\beta^\vee \in 2Q^\vee$, then $\alpha,\beta$ belong to the same copy of $G_2$. \end{lem} \pf Let $V$ denote the real vector space spanned by
$R$. Choose a Weyl-invariant scalar product $(\ ,\ )$ and use it to identify $V$ with its dual and regard $R^\vee$ as a root system in $V$.

Assume now that $\alpha^\vee+\beta^\vee \in 2Q^\vee$. Note that $\alpha^\vee$ and $\beta^\vee$ are orthogonal (but may not be strongly orthogonal elements of
$R^\vee$).

First we show that then $\alpha,\beta$ belong to the same irreducible piece of $R$. To that end, assume
that $R$ decomposes as $R=R_1 \sqcup R_2$ and $V$ decomposes accordingly as $V_1 \oplus V_2$. If $\alpha \in R_1$ and $\beta \in R_2$, then $\alpha^\vee \in V_1$
and $\beta^\vee \in V_2$. Then $\frac{1}{2}(\alpha^\vee + \beta^\vee) \in Q^\vee$ implies $\frac{1}{2}\alpha^\vee \in Q_1^\vee,\frac{1}{2}\beta^\vee \in Q_2^\vee$
(project orthogonally onto $V_1$ resp $V_2$). This however contradicts the above lemma, because $\frac{1}{2}\alpha^\vee$ has length strictly less then the shortest
elements in $R_1^\vee$.

Knowing that $\alpha,\beta$ lie in the same irreducible piece we can now assume wlog that $R$ is irreducible. Normalize $(\ ,\ )$ so that the short roots in $R$
have length $1$. We have the following cases
\begin{itemize}
\item All elements of $R$ have length $1$. Then all elements of $R^\vee$ have length $2$. The length of $\frac{1}{2}(\alpha^\vee+\beta^\vee)$
is $\sqrt{2}$, which by the above lemma is not a length of an element in $Q^\vee$.

\item $R$ contains elements of lengths $1$ and $\sqrt{2}$. Then $R^\vee$ contains elements of lengths $\sqrt{2}$ and $2$.
\begin{itemize}
\item If both $\alpha^\vee,\beta^\vee$ have length $\sqrt{2}$, then $\frac{1}{2}(\alpha^\vee+\beta^\vee)$ has length $1$, so is not in $Q^\vee$.
\item If $\alpha^\vee$ has length $\sqrt{2}$ and $\beta^\vee$ has length $2$, then $\frac{1}{2}(\alpha^\vee+\beta^\vee)$ has length $\frac{\sqrt{6}}{2}$, so again is not in $Q^\vee$.
\item If both $\alpha^\vee,\beta^\vee$ have length $2$, then $\frac{1}{2}(\alpha^\vee+\beta^\vee)$ has length $\sqrt{2}$ and thus could potentially be in $Q^\vee$. If it is,
then by the above lemma it is also in $R^\vee$, so $\frac{1}{2}(\alpha^\vee+\beta^\vee)^\vee$ must be an element of $R$. One immediately computes that $[\frac{1}{2}(\alpha^\vee+\beta^\vee)]^\vee = \alpha+\beta$, but the latter is not
an element of $R$ because $\alpha,\beta$ are strongly orthogonal.
\end{itemize}

\item $R$ has elements of lengths $1$ and $\sqrt{3}$. Then $R$ is $G_2$ and $R^\vee$ is also $G_2$. As one sees immediately, up to the action of its Weyl-group,
$G_2$ has a unique pair of orthogonal roots, which are then automatically strongly orthogonal and half their sum is also a root.
\end{itemize}
\qed

%---------------------------------------------------------------------------------
% Section
%=================================================================================

\section{Splitting invariants} \label{sec:splitting}
Recall that we have fixed a split semi-simple and simply-connected group $G$ over $\R$ and a splitting $(T_0,B_0,\{X_\alpha\})$ of it. Given a maximal torus
$T$, an element $h \in G$ s.t. $\tx{Int}(h)T_0=T$ and a-data $\{a_\beta\}$ for $R(T,G)$, Langlands and Shelstad construct in \cite[2.3]{LS1} a certain element
of $Z^1(\Gamma, T)$, whose image in $H^1(\Gamma,T)$ they call $\lambda(T)$ -- the "splitting invariant" of $T$. They show that this image is independent of the
choice of $h$. In this section we want to study this splitting invariant in such a way that enables us to see how it varies when the torus varies. It turns out
that a certain type of a-data is very well suited for this. This a-data is determined by a Borel $B \supset T$ as follows:
\[ \alpha_\beta = \begin{cases}i&,\beta>0 \wedge \sigma_T(\beta)<0 \\-i&,\beta<0 \wedge \sigma_T(\beta)>0 \\ \
                               1&,\beta>0 \wedge \sigma_T(\beta)>0 \\-1&,\beta<0 \wedge \sigma_T(\beta)<0 \end{cases} \]
where $\sigma_T$ denotes the Galois-action on $X^*(T)$ and $\beta>0$ means $\beta \in R(T,B)$. We will call this a-data $B$-a-data. It should not be confused
with Shelstad's terminology of \ti{based} a-data, which is also given by a Borel -- for based a-data, the positive imaginary roots are assigned $i$ while all
other positive roots are assigned $1$; for $B$-a-data, any positive root whose Galois-conjugate is negative is assigned $i$. Therefore a splitting invariant
computed using Shelstad's based a-data will in general be different from one computed using $B$-a-data. The precise difference is given by \cite[2.3.2]{LS1}.
It is however more important to note that according to \cite[Lemma 2.3.C]{LS1} this difference disappears once the splitting invariant has been paired with an
endoscopic character. Thus, as far applications to transfer factors are concerned, based a-data and $B$-a-data give the same result.

 In view of the reduction theorem which we will prove in section \ref{sec:splitting_A}, it will be helpful to consider not just the cohomology class, but
the actual cocycle constructed in \cite[2.3]{LS1}. We will denote this cocycle by $\lambda(T,B,h)$ to record its dependence on the $B$-a-data and the element
$h$, while the splitting $(T_0,B_0,\{X_\alpha\})$ is not present in the notation because it is assumed fixed. Since we are working over $\R$, we will identify
a 1-cocycle and its value at $\sigma \in \tx{Gal}(\C/\R)$, and hence we will view $\lambda(T,B,h)$ as an element of $T$. Given $h$, there is an obvious choice
for $B$, namely $\tx{Int}(h)B_0$. We will write $\lambda(T,h)$ for $\lambda(T,\tx{Int}(h)B_0,h)$. Note that in this notation, $T$ is clearly redundant, because
it equals $\tx{Int}(h)T_0$. However, we keep it so that the notation is close to that in \cite{LS1}. We would like to alert the reader of one potential
confusion -- while the cohomology class of $\lambda(T,B,h)$ is independent of the choice of $h$, that of $\lambda(T,h)$ is not, because in the latter $h$
influences not only the identification of $T_0$ with $T$ but also the choice of $B$-a-data for $T$.

%---------------------------------------------------------------------------------
% Subsection
%---------------------------------------------------------------------------------

\subsection{The splitting invariant for $T_\alpha$}

Recall from section \ref{sec:cayley} that for each $\alpha \in R(T_0,B_0)$ there is a canonical maximal torus $T_\alpha$ and a pair of isomorphisms
$T_0^{s_\alpha} \rightarrow T_\alpha$. To fix one of the two, fix $\mu \in \Omega$ s.t. $\mu^{-1}\alpha \in \Delta$. Then $\tx{Int}(g_{\mu,\alpha})$ is one of
the two isomorphisms $T_0^{s_\alpha} \rightarrow T_\alpha$. The goal of this section is to compute $\lambda(T_\alpha,B,g_{\mu,\alpha})$ for a given Borel $B
\supset T_\alpha$, and in particular $\lambda(T_\alpha,g_{\mu,\alpha})$. We will give a formula for the latter in purely root-theoretic terms.

%\begin{dfn} For $h \in G$, a Borel $B$ of $T:=\tx{Int}(h)T_0$, and $\alpha \in R$ put
%\[ \epsilon_{h,B}(\alpha) = \begin{cases}+1&,\alpha\circ\tx{Int}(h^{-1}) \in R(T,B)\\-1&,\tx{else}\end{cases} \]
%\end{dfn}

\begin{lem}\label{lem:lambda_one} With $g:=g_{\mu,\alpha}$ we have
\[ \lambda(T_\alpha,B,g) = \tx{Int}(g)\Big(\alpha^\vee(i\cdot a_{\alpha \circ \tx{Int}(g^{-1})}) \cdot s_\alpha(\sigma(\delta))\delta^{-1} \Big)\]
where
\[ \delta = \prod_{\substack{\beta>0\\\mu^{-1}\beta<0}} \beta^\vee( a_{\beta\circ \tx{Int}(g^{-1})} )^{-1} \]
and $\sigma$ denotes complex conjugation on $T_0$.
\end{lem}
\pf Put $u=n(\mu)$. We will first compute the cocycle $\lambda(T_\alpha,B,g u)$. The notation will be as in \cite[2.3]{LS1}. The pullback of the
$\Gamma$-action on $T_\alpha$ to $T_0$ via $gu$ differs from $\sigma$ by
\begin{eqnarray*}
\omega_{T_\alpha}(\sigma) &:=& \tx{Int}\left( (gu)^{-1} \sigma(gu) \right)\
= \tx{Int}\left( n(\mu)^{-1} g^{-1} \sigma(g) n(\mu) \right)\
= \mu^{-1}s_\alpha\mu\\
&=& s_{\mu^{-1}\alpha}
\end{eqnarray*}
Using that $\mu^{-1}\alpha$ is simple, we compute the three factors of $\tx{Int}(gu)^{-1}\lambda(T_\alpha,B,gu)$:
\begin{eqnarray*}
x(\sigma)&=&\prod_{\substack{\beta>0\\ \omega_{T_\alpha}(\sigma)\beta<0}} \beta^\vee( a_{\beta\circ\tx{Int}(gu)^{-1}} )\\
&=&(\mu^{-1}\alpha)^\vee(a_{\mu^{-1}\alpha\circ\tx{Int}(u^{-1}g^{-1})})\\
%&=&(\mu^{-1}\alpha)^\vee(a_{\alpha\circ\tx{Int}(g^{-1})})\\
&=&\tx{Int}(u^{-1})(\alpha^\vee(a_{\alpha\circ\tx{Int}(g^{-1})}))\\ \\
n(\omega_{T_\alpha}(\sigma))&=&\mat{0}{1}{-1}{0}_{X_{\mu^{-1}\alpha}} = \tx{Int}(u^{-1})\mat{0}{1}{-1}{0}_{X_{\mu|\mu^{-1}\alpha}}\\ \\
\sigma(gu)^{-1}(gu)&=&\tx{Int}(u^{-1})g^2 = \tx{Int}(u^{-1})\mat{0}{i}{i}{0}_{X_{\mu|\mu^{-1}\alpha}}
\end{eqnarray*}
Thus
\[ \lambda(T_\alpha,B,gu) = \tx{Int}(gu)\tx{Int}(u^{-1})\Big(\alpha^\vee(a_{\alpha\circ\tx{Int}(g^{-1})})\alpha^\vee(i)\Big)\]
From the proofs of \cite[2.3.A]{LS1} and \cite[2.3.B]{LS1} one sees that
\[ \lambda(T_\alpha,B,gu) = \tx{Int}(g)\Big(\delta\sigma_{T_\alpha}(\delta)^{-1}\Big)\cdot\lambda(T_\alpha,B,g) \]
where $\sigma_{T_\alpha}$ is the transport of the action of complex conjugation on $T_\alpha$ to $T_0$ via $g$. This action is $s_\alpha \rtimes \sigma$.
Notice that the term $\lambda^{-1}\sigma_T(\lambda)$ appearing in the proof of \cite[2.3.A]{LS1} is trivial since for us $u=n(\mu)$ and hence $\lambda=1$. The
claim now follows. \qed

Before we turn to the computation of $s_\alpha(\sigma(\delta))\delta^{-1}$ we will need to take a closer look at the following set.
\begin{dfn}
For $\alpha>0$ put $R^+_\alpha=\Big\{ \beta \in R|\ \beta>0 \wedge s_\alpha(\beta)<0 \Big\}$
\end{dfn}

\begin{lem} \label{lem:ralphasplit} Let $\alpha>0$ and $\mu \in \Omega$ be s.t. $\mu^{-1}\alpha \in \Delta$. Then the sets
\begin{eqnarray*}
&&\Big\{\beta \in R|\ \beta>0 \wedge s_\alpha(\beta)<0 \wedge \mu^{-1}\beta <0 \Big\}\\
\tx{and}&&\\
&&\Big\{\beta \in R|\ \beta>0 \wedge s_\alpha(\beta)<0 \wedge \mu^{-1}\beta >0 \wedge \beta\neq\alpha\Big\}
\end{eqnarray*}
are disjoint and their union is $R^+_\alpha - \{\alpha\}$. The map
\[ \beta \mapsto -s_\alpha(\beta) \]
is an involution on $R^+_\alpha-\{\alpha\}$ which interchanges the above two sets.
\end{lem}
\pf Every $\beta$ in the first set satisfies $\beta \neq\alpha$ because $\mu^{-1}\alpha$ is positive. Hence the first set lies in $R^+_\alpha-\{\alpha\}$ and
clearly the second also does. The fact that the two are disjoint and cover $R^+_\alpha-\{\alpha\}$ is obvious. Now to the bijection. Let $\beta$ be an
element in the first set, and consider $\tilde\beta=-s_\alpha(\beta)$.
We have
\begin{eqnarray*}
\beta\neq\alpha&\Rightarrow&\tilde\beta\neq\alpha\\
s_\alpha(\beta)<0&\Rightarrow&\tilde\beta>0\\
\beta>0&\Rightarrow&s_\alpha\tilde\beta=-\beta<0\\
\mu^{-1}\beta<0&\Rightarrow&\mu^{-1}\tilde\beta=\mu^{-1}s_\alpha(-\beta)=s_{\mu^{-1}\alpha}(-\mu^{-1}\beta)>0
\end{eqnarray*}
where the last inequality holds because $\mu^{-1}\alpha$ is simple and
\[ \beta>0\Rightarrow\beta\neq-\alpha\Rightarrow-\mu^{-1}\beta \neq \mu^{-1}\alpha \]
\qed

\rmk A similar observation appears in \cite[\S4.3]{LS2}.

\begin{lem} \label{lem:delta} We have
\[ s_\alpha(\sigma(\delta))\delta^{-1} = \prod_{\substack{\beta\in R^+_\alpha\\\mu^{-1}\beta<0}} \Big[ \beta^\vee\big(a_{\beta \circ \tx{Int}(g^{-1})}\big)\ s_\alpha\beta^\vee\big(a_{s_\alpha\beta \circ \tx{Int}(g^{-1})}\big)^{-1}\Big] \]
\end{lem}
\pf According to the proof of part (a) of \cite[2.3.B]{LS1} the contributions to $\delta s_\alpha(\sigma(\delta))^{-1}$ are as follows:
\begin{enumerate}
\item[(1)] $\{\beta|\ \beta>0\wedge\mu^{-1}\beta<0 \wedge s_\alpha\beta<0 \} : \beta^\vee(a_{\beta\circ\tx{Int}(g^{-1})})^{-1}$
\item[(2)] $\{\beta|\ \beta<0\wedge\mu^{-1}\beta<0 \wedge \mu^{-1}s_\alpha\beta<0 \wedge s_\alpha\beta>0 \} : \beta^\vee(a_{\beta\circ\tx{Int}(g^{-1})})$
\item[(3)] $\{\beta|\ \beta>0\wedge\mu^{-1}\beta<0 \wedge s_\alpha\beta>0 \wedge \mu^{-1}s_\alpha\beta>0 \} :
    \beta^\vee(a_{\beta\circ\tx{Int}(g^{-1})})^{-1}$
\item[(4)] $\{\beta|\ \mu^{-1}\beta>0 \wedge s_\alpha\beta>0 \wedge \mu^{-1}s_\alpha\beta <0 \} : \beta^\vee(a_{\beta\circ\tx{Int}(g^{-1})})$
\end{enumerate}
We will use $\mu^{-1}s_\alpha(\beta) = s_{\mu^{-1}\alpha}(\mu^{-1}\beta)$ and the fact that $\mu^{-1}\alpha$ is simple to show that the last two sets are empty.
In set $(3)$, the conditions $\mu^{-1}\beta<0$ and $\mu^{-1}s_\alpha\beta>0$ imply $\mu^{-1}\beta=-\mu^{-1}\alpha$, i.e. $\beta=-\alpha$, which contradicts $\beta>0$.
In set $(4)$, the conditions $\mu^{-1}\beta>0$ and $\mu^{-1}s_\alpha\beta<0$ imply $\beta=\alpha$. Since $\alpha>0$ this contradicts $s_\alpha\beta>0$.

Next we claim $(2)=s_\alpha((1))$. We have
\[ \mu^{-1}\beta<0 \wedge \mu^{-1}s_\alpha\beta<0 \Leftrightarrow \mu^{-1}\beta<0\wedge\mu^{-1}\beta\neq-\mu^{-1}\alpha \]
from which we get
\[ (2) = \{ -\beta|\ \beta>0 \wedge s_\alpha\beta<0 \wedge \mu^{-1}\beta>0 \wedge \beta \neq \alpha \} \]
Now $(2)=s_\alpha((1))$ follows from Lemma \ref{lem:ralphasplit}.

From these considerations it follows that
\begin{eqnarray*}
\delta s_\alpha(\sigma(\delta))^{-1}&=&\prod_{\beta\in(1)}\beta^\vee(a_{\beta\circ\tx{Int}(g^{-1})})^{-1} \prod_{\beta\in(2)}\beta^\vee(a_{\beta\circ\tx{Int}(g^{-1})})\\
&=&\prod_{\beta \in (1)} \Big[ \beta^\vee\big(a_{\beta \circ \tx{Int}(g^{-1})}\big)^{-1}\ s_\alpha\beta^\vee\big(a_{s_\alpha\beta \circ \tx{Int}(g^{-1})}\big)\Big]
\end{eqnarray*}
\qed

Let us recall our notation: $\alpha \in R$ is any positive root, $\mu \in \Omega$ is s.t. $\mu^{-1}\alpha \in \Delta$, and $g=g_{\mu,\alpha}$ is the Cayley-transform
corresponding to $X_{\mu|\mu^{-1}\alpha}$.

From Lemmas \ref{lem:lambda_one} and \ref{lem:delta} we immediately get
\begin{cor}
\[ \lambda(T_\alpha,B,g) = \tx{Int}(g)\left( \alpha^\vee(ia_{\alpha \circ\tx{Int}(g^{-1})})\cdot\prod_{\substack{\beta\in R^+_\alpha\\\mu^{-1}\beta<0}}\
   \Big[ \beta^\vee\big(a_{\beta \circ \tx{Int}(g^{-1})}\big)\ s_\alpha\beta^\vee\big(a_{s_\alpha\beta \circ \tx{Int}(g^{-1})}\big)^{-1}\Big] \right) \]
\end{cor}

In the case $B=\tx{Int}(g)B_0$ this formula becomes simpler.

\begin{cor} \label{cor:lambda}
\[ \lambda(T_\alpha,g_{\mu,\alpha}) = \tx{Int}(g_{\mu,\alpha})\left( \alpha^\vee(-1)\cdot\prod_{\substack{\beta\in R^+_\alpha\\\mu^{-1}\beta<0}}(\beta^\vee\cdot s_\alpha\beta^\vee)(i)\right) \]
\end{cor}

\begin{dfn} Put
\[ \rho(\mu,\alpha) := \alpha^\vee(-1)\cdot\prod_{\substack{\beta\in R^+_\alpha\\\mu^{-1}\beta<0}}(\beta^\vee\cdot s_\alpha\beta^\vee)(i) \quad \in T_0\]
\end{dfn}
By Corollary \ref{cor:lambda} and the work of \cite[2.3]{LS1} we know that $\rho(\mu,\alpha) \in Z^1(\Gamma,T_0^{s_\alpha})$ and
$\tx{Int}(g_{\mu,\alpha})\rho(\mu,\alpha) = \lambda(T_\alpha,g_{\mu,\alpha})$.

\begin{pro}\ \\[-20pt] \label{pro:changemu}
\begin{enumerate}
\item $\rho(\mu,\alpha)=\prod\limits_{\beta \in R_\alpha^+} \beta^\vee(i)n(s_\alpha)g_{\mu,\alpha}^2$.
\item $s_\alpha\rho(\mu,\alpha)=\rho(\mu,\alpha)$, $\sigma(\rho(\mu,\alpha))=\rho(\mu,\alpha)^{-1}$.
\item The image of $\rho(\mu,\alpha)$ under the two canonical isomorphisms $T_0^{s_\alpha} \rightarrow T_\alpha$ is the same
\item If $\mu' \in \Omega$ is another Weyl-element s.t. $(\mu')^{-1}\alpha \in \Delta$, then
\[ \rho(\mu',\alpha) = \alpha^\vee\Big(\epsilon(\mu',\alpha,\mu)\Big)\rho(\mu,\alpha) \]
In particular
\[ \lambda(T_\alpha,g_{\mu',\alpha}) = \lambda(T_\alpha,g_{\mu,\alpha}) \cdot \tx{Int}(g_{\mu,\alpha})\Big[ \alpha^\vee(\epsilon(\mu',\alpha,\mu))\Big]  \]
\end{enumerate}
\end{pro}
\pf The first point follows from Corollary \ref{cor:lambda}, because the right hand side is by construction
$\tx{Int}(g_{\mu,\alpha})\lambda(T_\alpha,g_{\mu,\alpha})$. The second point is evident from the structure of $\rho$. The third point is now clear because as
remarked in section \ref{sec:cayley} the two canonical isomorphisms differ by precomposition with $s_\alpha$.

For the last point,
\[ \rho(\mu',\alpha) = \prod_{\beta\in R^+_\alpha} \beta^\vee(i) n(s_\alpha) g_{\mu',\alpha}^2 \]
If $\epsilon(\mu',\alpha,\mu)=+1$ then $g_{\mu',\alpha}=g_{\mu,\alpha}$ and the statement is clear. Assume now that $\epsilon(\mu',\alpha,\mu)=-1$. Then
$g_{\mu',\alpha}=g_{\mu,\alpha}^{-1}$. We see
\[ \rho(\mu',\alpha)=\prod_{\beta\in R^+_\alpha} \beta^\vee(i) n(s_\alpha) g_{\mu,\alpha}^2g_{\mu,\alpha}^{-4} \]
But $g_{\mu,\alpha}^{-4}=\alpha^\vee(-1)$, hence the claim.  \qed

%---------------------------------------------------------------------------------
% Subsection
%---------------------------------------------------------------------------------

\subsection{The splitting invariant for $T_A$} \label{sec:splitting_A}

\begin{fct} Let $A$ be a SOS in $R$. Consider the set of automorphisms of $G$ given by
\[ \left\{ \tx{Int}(g) \Bigg|\ g = \prod_{\alpha \in A} g_{\mu_\alpha,\alpha}, \mu_\alpha \in \Omega, \mu_\alpha^{-1}\alpha \in \Delta \right\} \]
The image of $T_0$ under any element of that set is the same. Call it $T_A$. Then any element of that set induces an isomorphism of real tori
\[ T_0^{S_A} \rightarrow T_A \]
\end{fct}

\pf Let $\tx{Int}(g_1),\tx{Int}(g_2)$ be elements of the above set, with
\[ g_i = \prod_{\alpha \in A} g_{\mu^i_\alpha,\alpha} \]
and let $A' \subset A$ be the subset of those $\alpha$ s.t. $g_{\mu^1_\alpha,\alpha} \neq g_{\mu^2_\alpha,\alpha}$. For those $\alpha$ we have then
$g_{\mu^1_\alpha,\alpha}=g_{\mu^2_\alpha,\alpha}^{-1}$, hence $\tx{Int}(g_1g_2^{-1})|_{T_0} = \prod_{\alpha \in A'}g_{\mu^1_\alpha,\alpha}^2|_{T_0}=S_{A'}$
which normalizes $T_0$. This shows that the images of $T_0$ under these two automorphisms are the same. Moreover, the transport of the $\Gamma$-action on $T_A$
to $T_0$ via $\tx{Int}(g_1^{-1})$ differs from $\sigma$ by $\tx{Int}(\sigma(g_1^{-1})g_1)|_{T_0}=\tx{Int}(g_1^2)|_{T_0}=S_A$. \qed

\begin{dfn} \label{dfn:canonicalset} For a SOS $A \subset R$, we will call the set
\[ \left\{ \tx{Int}(g)|_{T_0} \Bigg|\ g = \prod_{\alpha \in A} g_{\mu_\alpha,\alpha}, \mu_\alpha \in \Omega, \mu_\alpha^{-1}\alpha \in \Delta \right\} \]
the canonical set of isomorphisms $T_0^{S_A} \rightarrow T_A$. More generally, if $A' \subset A$, we will call the set
\[ \left\{ \tx{Int}(g)|_{T_{A'}} \Bigg|\ g = \prod_{\alpha \in A \setminus A'} g_{\mu_\alpha,\alpha}, \mu_\alpha \in \Omega, \mu_\alpha^{-1}\alpha \in \Delta \right\} \]
the canonical set of isomorphisms $T_{A'}^{S_{A \setminus A'}} \rightarrow T_A$.
\end{dfn}

\begin{fct} Any maximal torus $T \subset G$ is $G(\R)$-conjugate to one of the $T_A$. \end{fct}
\pf Choose $g \in G$ s.t. $\tx{Int}(g)T_0=T$. The transport of the $\Gamma$-action on $T$ to $T_0$ via $\tx{Int}(g^{-1})$ differs from $\sigma$ by an element
of $Z^1(\Gamma,\Omega)=\tx{Hom}(\Gamma,\Omega)$ and this element sends complex conjugation to an element of $\Omega$ of order 2. By \cite[Ch.VI.Ex
\S1(15)]{Bou} there exists a SOS $A$ s.t. this element equals $S_A$. If $\tx{Int}(g_A)$ is one of the canonical isomorphisms $T_0^{S_A} \rightarrow T_A$, then
$\tx{Int}(g_Ag^{-1}) : T \rightarrow T_A$ is an isomorphism of real tori. By \cite[Thm. 2.1]{S1} there exists $g' \in G(\R)$ s.t. $\tx{Int}(g')T=T_A$.\qed

If we conjugate $A$ by $\Omega$ to an $A'$, then the tori $T_A$ and $T_{A'}$ are also conjugate by $G(\R)$. Thus we may fix representatives $A_1,...,A_k$ for
the $\Omega$-orbits of MSOS in $R$ and study the tori $T_A$ for $A$ inside one of the $A_i$. We assume that the fixed splitting $(T_0,B_0,\{X_\alpha\})$ is
compatible with the choice of representatives in the following sense
\begin{itemize}
\item $\#\#(R,>,A_i)$ holds for all $A_i$.
\item If $\alpha,\beta \in A_i$ lie in the same $G_2$-factor then one of them is simple
\end{itemize}
This can always be arranged, as Lemmas \ref{lem:msosreps}, \ref{lem:posroots}, Fact \ref{fct:g2sos} and Proposition \ref{pro:strong} show. Notice that this
condition does not reduce generality -- it is only a condition on $B_0$, but all Borels containing $T_0$ are conjugate under $N_{T_0}(\R)$ and thus by
\cite[2.3.1]{LS1} the splitting invariants are independent of the choice of $B_0$.

%Recall from section \ref{sec:cayley} that once a splitting is fixed, we get for each root $\alpha$ a canonical pair of elements $g_\alpha,g_\alpha^{-1} \in G$
%with the property that $g_\alpha^2,g_\alpha^{-2} \in N(T_0)$ and $\tx{Int}(g_\alpha^2)|_{T_0} = \tx{Int}(g_\alpha^{-2})|_{T_0} = s_\alpha$. Each of these
%elements is given as $g_{\mu,\alpha}$ for some $\mu \in \Omega$ s.t. $\mu^{-1}\alpha \in \Delta$. Moreover, if $\alpha,\beta$ are strongly-orthogonal then
%$g_\alpha,g_\beta$ commute by Fact \ref{fct:strongcommute}.

\begin{lem} \label{lem:commute} If $A',A''$ are disjoint subsets of some $A_i$ then
\[ n(S_{A'})n(S_{A''}) = n(S_{A' \cup A''})\]
In particular, $n(S_{A'})$ and $n(S_{A''})$ commute.
\end{lem}
\pf This follows immediately from \cite[Lemma 2.1.A]{LS1}, because by $\#\#$ the set
\[ \{\beta \in R:\ \beta>0 \wedge S_{A'}(\beta)<0 \wedge S_{A'}S_{A''}(\beta)>0 \} \]
is empty. \qed

%\pf The above set equals
%\[ \{ \beta \in R^+_{A'}:\ S_{A'}S_{A''}(\beta)>0 \} \]
%Let $\beta \in R^+_{A'}$. Then $\beta' := -s_{A'}(\beta) \in R^+_{A'}$ and since $S_{A''}$ preserves $R^+_{A'}$ we have
%\[ S_{A''}S_{A'}(\beta) = -S_{A''}\beta' < 0 \]
%is empty. \qed

%\pf
%\begin{eqnarray*}
%&&\\
%&=&
%\end{eqnarray*}

\begin{pro}\label{pro:reflectrho}
Let $\alpha,\gamma$ be distinct elements of one of the $A_i$, and $\mu \in \Omega$ be s.t. $\mu^{-1}\alpha \in \Delta$. Then $\rho(\mu,\alpha)$ is fixed by $s_\gamma$.
\end{pro}

\pf We first show\\
\ul{Claim}:
\[ s_\gamma\rho(\mu,\alpha) = \rho(\mu,\alpha)\alpha^\vee(\epsilon(s_\gamma\mu,\alpha,\mu)) \]
\pf We have
\[ s_\gamma(\rho(\mu,\alpha)) = s_\gamma\left(\alpha^\vee(-1)\prod_{\substack{\beta\in R^+_\alpha\\ \mu^{-1}\beta < 0}} \beta^\vee(i)s_\alpha\beta^\vee(i) \right) \]
Now $s_\gamma$ preserves $\alpha^\vee$, commutes with $s_\alpha$, and by Proposition \ref{pro:decomp} also preserves the set $R^+_\alpha$,
hence the last expression equals
\begin{eqnarray*}
\alpha^\vee(-1) \prod_{\substack{\beta\in R^+_\alpha\\ \mu^{-1}\beta < 0}} s_\gamma\beta^\vee(i)s_\alpha s_\gamma\beta^\vee(i)\
&=&\alpha^\vee(-1) \prod_{\substack{s_\gamma\beta\in R^+_\alpha\\ \mu^{-1}s_\gamma\beta < 0}} \beta^\vee(i)s_\alpha\beta^\vee(i)\\
&=&\alpha^\vee(-1) \prod_{\substack{\beta\in R^+_\alpha\\ \mu^{-1}s_\gamma\beta < 0}} \beta^\vee(i)s_\alpha\beta^\vee(i)\\
&=&\rho(s_\gamma\mu,\alpha)\\
&=&\rho(\mu,\alpha) \cdot \alpha^\vee\Big(\epsilon(s_\gamma\mu,\alpha,\mu)\Big)
\end{eqnarray*}
the last equality coming from Proposition \ref{pro:changemu}. \qed(claim)

We want to show $\alpha^\vee\Big(\epsilon(s_\gamma\mu,\alpha,\mu)\Big)=1$. Choose $\nu\in\Omega$ s.t. $\nu^{-1}\gamma \in \Delta$. We will derive and compare
two expressions for
\begin{equation}\label{eq:one} \prod_{\beta \in R^+_{\{\alpha,\gamma\}}}\beta^\vee(i) n(s_\alpha s_\gamma)g_{\mu,\alpha}^2g_{\nu,\gamma}^2 \end{equation}
By Corollary \ref{cor:decomp} we have
\[\prod_{\beta \in R^+_{\{\alpha,\gamma\}}}\beta^\vee(i) = \prod_{\beta \in R^+_\alpha} \beta^\vee(i) \prod_{\beta \in R^+_\gamma} \beta^\vee(i) \]
By Proposition \ref{pro:decomp} $s_\alpha$ is a permutation of the set $R^+_\gamma$, hence
\[ n(s_\alpha)\prod_{\beta \in R^+_\gamma} \beta^\vee(i) n(s_\alpha)^{-1}  = \prod_{\beta \in R^+_\gamma} \beta^\vee(i) \]
By Lemma \ref{lem:commute}, the elements $n(s_\alpha)$ and $n(s_\gamma)$ of $N(T_0)$ commute. Moreover, by Fact \ref{fct:strongcommute} the elements
$g_{\mu,\alpha}^2$ and $g_{\nu,\gamma}^2$ commute. Thus we get on the one hand
\begin{eqnarray*}
\eqref{eq:one}&=&\prod_{\beta \in R^+_\gamma} \beta^\vee(i) n(s_\gamma) \prod_{\beta \in R^+_\alpha} \beta^\vee(i)n(s_\alpha) g_{\mu,\alpha}^2g_{\nu,\gamma}^2\\
&=&\rho(\nu,\gamma) \tx{Int}(g_{\nu,\gamma}^{-2})\Big[\rho(\mu,\alpha)\Big]\\
&=&\rho(\nu,\gamma)\rho(\mu,\alpha)\alpha^\vee(\epsilon(s_\gamma\mu,\alpha,\mu))
\end{eqnarray*}
where the last equality follows from above claim. Analogously, we obtain on the other hand
\begin{eqnarray*}
\eqref{eq:one}&=&\prod_{\beta \in R^+_\alpha} \beta^\vee(i) n(s_\alpha) \prod_{\beta \in R^+_\gamma} \beta^\vee(i)n(s_\gamma) g_{\nu,\gamma}^2g_{\mu,\alpha}^2\\
&=&\rho(\mu,\alpha) \tx{Int}(g_{\mu,\alpha}^{-2})\Big[\rho(\nu,\gamma)\Big]\\
&=&\rho(\mu,\alpha)\rho(\nu,\gamma)\gamma^\vee(\epsilon(s_\alpha\nu,\gamma,\nu))
\end{eqnarray*}

We conclude that
\[ \alpha^\vee(\epsilon(s_\gamma\mu,\alpha,\mu)) = \gamma^\vee(\epsilon(s_\alpha\nu,\gamma,\nu)) \]
We claim that both sides of this equality are trivial. Assume by way of contradiction that this is not the case. Then we have
\begin{eqnarray*}
\alpha^\vee(-1)=\gamma^\vee(-1)&\Leftrightarrow&1 = (-1)^{(\alpha^\vee-\gamma^\vee)} = (-1)^{\alpha^\vee+\gamma^\vee} \in \C^\times \otimes X_*(T_0) = \C^\times \otimes Q^\vee\\
&\Leftrightarrow&\alpha^\vee+\gamma^\vee \in 2Q^\vee
\end{eqnarray*}
where $Q^\vee$ is the coroot-lattice of $T_0$, which coincides with $X_*(T_0)$ since $G$ is simply-connected. By Lemma \ref{lem:strongco} $\alpha,\gamma$ must lie
in the same $G_2$-factor of $R$. In this case $\{\alpha,\gamma\}$ is a MSOS for that $G_2$-factor and by our assumption from the beginning of this section one of
$\alpha,\gamma$ must be simple. Say wlog $\alpha$ is simple. By Proposition \ref{pro:changemu}
\[\rho(\mu,\alpha)=\alpha^\vee(\epsilon(\mu,\alpha,1))\rho(1,\alpha) = \alpha^\vee(-\epsilon(\mu,\alpha,1)) \]
which is clearly fixed by $s_\gamma$, thus we see
\[ 1 = s_\gamma(\rho(\mu,\alpha))\rho(\mu,\alpha)^{-1} = \alpha^\vee(\epsilon(s_\gamma\mu,\alpha,\mu)) = \gamma^\vee(\epsilon(s_\alpha\nu,\gamma,\nu)) \]

\qed

\begin{cor} Let $A$ be a subset of some $A_i$, $\alpha \in A$ and $\mu \in \Omega$ s.t. $\mu^{-1}\alpha \in \Delta$. Then $\rho(\mu,\alpha) \in
Z^1(\Gamma,T_0^{S_A})$ and its image in $T_A$ under any canonical isomorphism $T_0^{S_A} \rightarrow T_A$ is the same. \end{cor}
\pf By Propositions \ref{pro:changemu} and \ref{pro:reflectrho} $\rho=\rho(\mu,\alpha)$ is fixed by $s_\gamma$ for any $\gamma \in A$. The first statement now
follows from $\rho S_A\sigma(\rho) = \rho \sigma(\rho)=1$ showing $\rho \in Z^1(\Gamma,T_0^{S_A})$. The second holds because any two canonical isomorphisms
$T_0^{S_A} \rightarrow T_A$ differ by precomposition with $S_{A'}$ for some $A' \subset A$\qed

\begin{cor}Let $A$ be a subset of some $A_i$. Choose, for each $\alpha\in A$, $\mu_\alpha \in \Omega$ s.t. $\mu_\alpha^{-1}\alpha \in \Delta$. Put
\[ \rho(\{\mu_\alpha\}_{\alpha \in A},A) = \prod_{\alpha \in A} \rho(\mu_\alpha,\alpha) \]
Then
\begin{enumerate}
\item $\rho(\{\mu_\alpha\}_{\alpha \in A},A)$ is fixed by $s_\gamma$ for all $\gamma \in A$ (even all $\gamma \in A_i$)
\item The image of $\rho(\{\mu_\alpha\}_{\alpha \in A},A)$ under any of the canonical isomorphisms $T_0^{S_A} \rightarrow T_A$ is the same.
\end{enumerate}
\end{cor}
\pf Clear by the preceding Corollary.

\begin{pro} \label{pro:lambdadecomp}
Let $A$ be a subset of some $A_i$. For each $\alpha\in A$ choose $\mu_\alpha \in \Omega$ s.t. $\mu_\alpha^{-1}\alpha \in \Delta$ and put $g_A = \prod_{\alpha
\in A}g_{\mu_\alpha,\alpha}$. Then $\lambda(T_A,g_A)$ is the common image of $\rho(\{\mu_\alpha\}_{\alpha \in A},A)$ under the canonical isomorphisms
$T_0^{S_A} \rightarrow T_A$. In particular
\[ \lambda(T_A,g_A) = \prod_{\alpha \in A} \tx{Int}(g_{A-\{\alpha\}})\lambda(T_\alpha,g_{\mu_\alpha,\alpha}) \]
is a decomposition of $\lambda(T_A,g_A)$ as a product of elements of $Z^1(\Gamma,T_A)$.
\end{pro}
\pf The factors of the cocycle $\tx{Int}(g_A^{-1})\lambda(T_A,g_A) \in Z^1(\Gamma,T_0^{S_A})$ associated to these choices are as follows:
\[ x(\sigma_T)=\prod_{\beta\in R^+_A}\beta^\vee(i)=\prod_{\alpha\in A}\prod_{\beta \in R^+_\alpha}\beta^\vee(i) \]
where the second equality is due to Corollary \ref{cor:decomp},
\[ n(\omega_T(\sigma)) = n(S_A) = \prod_{\alpha \in A}n(s_\alpha) \]
where the second equality is due to Lemma \ref{lem:commute}, and
\[ \sigma(g_A)^{-1}g_A = \prod_{\alpha \in A}\sigma(g_\alpha)^{-1}g_\alpha = \prod_{\alpha \in A}g_\alpha^2 \]
Their product, which equals $\tx{Int}(g_A^{-1})\lambda(T_A,g_A)$, is thus
\[ x(\sigma_T)n(\omega_T(\sigma))\sigma(g_A)^{-1}g_A = \prod_{\alpha\in A}\prod_{\beta \in R^+_\alpha}\beta^\vee(i) \prod_{\alpha \in A}n(s_\alpha) \prod_{\alpha \in A}g_\alpha^2 \]
Just as in the proof of Proposition \ref{pro:reflectrho} we can rewrite this product as
\[ \prod_{\alpha\in A}\left[ \prod_{\beta \in R^+_\alpha}\beta^\vee(i) n(s_\alpha)\right] \prod_{\alpha \in A}g_\alpha^2 \]
Now we induct on the size of $A$, with $|A|=1$ being clear. Choose $\alpha_1 \in A$. Then
\begin{eqnarray*}
&\ &\prod_{\alpha \in A}\left[ \prod_{\beta \in R^+_{\alpha}}\beta^\vee(i) n(s_{\alpha})\right] \prod_{\alpha \in A}g_{\alpha}^2\\
&=&\prod_{\alpha \in A\setminus\{\alpha_1\}}\left[ \prod_{\beta \in R^+_{\alpha}}\beta^\vee(i) n(s_{\alpha})\right]  \left\{ \prod_{\beta \in R^+_{\alpha_1}}\beta^\vee(i) n(s_{\alpha_1}) g_{\alpha_1}^2 \right\} \prod_{\alpha \in A\setminus\{\alpha_1\}}g_{\alpha}^2\\
&=&\prod_{\alpha \in A\setminus\{\alpha_1\}}\left[ \prod_{\beta \in R^+_{\alpha}}\beta^\vee(i) n(s_{\alpha})\right] \rho(\mu_{\alpha_1},\alpha_1) \prod_{\alpha \in A\setminus\{\alpha_1\}}g_{\alpha}^2 \\
&=&\prod_{\alpha \in A\setminus\{\alpha_1\}}\left[ \prod_{\beta \in R^+_{\alpha}}\beta^\vee(i) n(s_{\alpha})\right] \prod_{\alpha \in A\setminus\{\alpha_1\}}g_{\alpha}^2\ \cdot\ \left(\prod_{\alpha \in A\setminus\{\alpha_1\}}^k s_{\alpha}\right)(\rho(\mu_{\alpha_1},\alpha_1))\\
&=&\prod_{\alpha \in A\setminus\{\alpha_1\}} \rho(\mu_{\alpha},\alpha)\ \cdot\ \rho(\mu_{\alpha_1},\alpha_1)\\
\end{eqnarray*}
where the last equality follows from Proposition \ref{pro:reflectrho} and the inductive hypothesis. This shows that
$\tx{Int}(g_A)\rho(\{\mu_\alpha\}_{\alpha\in A},A)=\lambda(T_A,g_A)$ and the result follows.\qed

\section{Explicit computations} \label{sec:explicit}
In this section we are going to use the classification of MSOS given in \cite{AK} to explicitly compute $\lambda(T_A,g_A)$ for the split simply-connected
semi-simple groups associated to the classical irreducible root systems. By Propositions \ref{pro:changemu} and \ref{pro:lambdadecomp} it is enough to compute
the cocycles $\rho(\mu,\alpha)$ for each $\alpha \in A_i$, and some $\mu \in \Omega$ with $\mu^{-1}\alpha \in \Delta$, where $A_1,...,A_k$ is a set of
representatives for the $\Omega$-classes of MSOS. We will use the notation from \cite{AK}, which is also the notation used in the Plates of \cite[Ch.VI]{Bou}.
There is only one cosmetic difference -- in \cite{Bou} the standard basis of $\R^k$ is denoted by $(\epsilon_i)$, in \cite{AK} by $(\lambda_i)$, and we are
using $(e_i)$. The dual basis will be denoted by $(e_i^*)$. One checks easily in each case that the choices of positive roots given in the Plates of
\cite[Ch.VI]{Bou} and the MSOS given in \cite{AK} satisfy condition $\#$ of section \ref{sec:strong}.

\subsection{Case $A_n$}
There is only one $\Omega$-equivalence class of MSOS, and the representative given in \cite{AK} is
\[ A_1 = \{ e_{2i-1}-e_{2i}|\ 1\leq i \leq [(n+1)/2] \} \]
All elements of this MSOS are simple roots and for each of them we can choose $\mu=1$. Then for any $\alpha \in A_1$ we have
\[ \boxed{\rho(1,\alpha) = \alpha^\vee(-1)} \]

\subsection{Case $B_n$}
If $n=2k+1$ then there is a unique equivalence class of MSOS, represented by
\[ A_1 = \{ e_{2i-1} \pm e_{2i}|\ 1 \leq i \leq k \} \cup \{ e_n \} \]
If $n=2k$ then there are two equivalence classes of MSOS, represented by
\begin{eqnarray*}
A_1&=&\{ e_{2i-1} \pm e_{2i} |\ 1 \leq i \leq k-1 \} \cup \{ e_n \}\\
A_2&=&\{ e_{2i-1}\pm e_{2i}|\ 1 \leq i \leq k \}
\end{eqnarray*}

If $\alpha = e_{2i-1}-e_{2i}$ or $\alpha = e_n$ then $\alpha$ is simple, we can choose $\mu=1$ and have
\[\boxed{\rho(1,\alpha)=\alpha^\vee(-1)}\]

If $\alpha = e_{2i-1}+e_{2i}$ then we take $\mu=s_{e_{2i}}$ and have $\mu^{-1}\alpha = e_{2i-1}-e_{2i} \in \Delta$. To compute $\rho(\mu,\alpha)$ we first observe
\begin{eqnarray*}
\{ \beta \in R|\ \beta>0 \wedge \mu^{-1}\beta <0 \wedge s_\alpha\beta<0 \}&=&\{ \beta \in R|\ \beta>0 \wedge \mu^{-1}\beta <0 \}\\
&=&\{e_{2i}\} \cup \{ e_{2i} \pm e_j|\ 2i<j \}
\end{eqnarray*}
hence
\begin{eqnarray*}
&&\sum_{\substack{\beta \in R^+_\alpha\\\mu^{-1}\beta<0}} (\beta^\vee+s_\alpha\beta^\vee)\\
&=&2e^*_{2i}+s_\alpha(2e^*_{2i}) + \sum_{j=2i+1}^n(e^*_{2i}-e^*_j+s_\alpha(e^*_{2i}-e^*_j) + e^*_{2i}+e^*_j+s_\alpha(e^*_{2i}+e^*_j))\\
&=&-2(e^*_{2i-1}-e^*_{2i})-2(e^*_{2i-1}-e^*_{2i})(n-2i)\\
&=&-2(n+1-2i)(e^*_{2i-1}-e^*_{2i})\\
&=&-2(n+1-2i)\alpha^\vee
\end{eqnarray*}
We get
\[ \boxed{\rho(\mu,\alpha) = \alpha^\vee( (-1)^n )} \]

\subsection{Case $C_n$}

The root system family $C_n$ is the only family for which the number of equivalence classes of MSOS grows when $n$ grows. Representatives for the equivalence
classes of MSOS are given by
\[ A_s = \{e_{2i-1}-e_{2i}|\ 1\leq i\leq s\} \cup \{2e_i|\ 2s+1 \leq i \leq n \} \qquad\qquad 0 \leq s \leq [n/2] \]
If $\alpha = e_{2i-1}-e_{2i}$ then $\alpha$ is simple and
\[ \boxed{\rho(1,\alpha)=\alpha^\vee(-1)} \]

If $\alpha = 2e_i$ and we take $\mu=s_{e_i-e_n}$ we have $\mu^{-1}\alpha = 2e_n \in \Delta$. Again we first observe
\begin{eqnarray*}
\{ \beta \in R|\ \beta>0 \wedge \mu^{-1}\beta <0 \wedge s_\alpha\beta<0 \}&=&\{ \beta \in R|\ \beta>0 \wedge \mu^{-1}\beta <0 \}\\
&=& \{ e_i-e_j |\ i<j \}
\end{eqnarray*}
hence
\[ \sum_{\substack{\beta \in R^+_\alpha\\\mu^{-1}\beta<0}} (\beta^\vee+s_\alpha\beta^\vee)=\sum_{j=i+1}^n(e^*_i-e^*_j)+(-e^*_i-e^*_j)\
=-2\sum_{j=i+1}^ne^*_j \]
We get
\[ \boxed{\rho(\mu,\alpha) = \prod_{j=i}^n e^*_i(-1)} \]

\subsection{Case $D_n$}

There is a unique equivalence class of MSOS represented by
\[ A_1 = \{ e_{2i-1} \pm e_{2i}|\ 1 \leq i \leq [n/2] \} \]
If $\alpha=e_{2i-1}-e_{2i}$ or $\alpha=e_{n-1}+e_n$ then $\alpha$ is simple and
\[ \boxed{\rho(1,\alpha)=\alpha^\vee(-1)} \]

If $\alpha=e_{2i-1}+e_{2i}$ with $2i \neq n$ then we take $\mu=s_{e_{2i-1}-e_{n-1}}\circ s_{e_{2i}-e_{n}}$ and have $\mu^{-1}\alpha = e_{n-1}-e_n \in \Delta$.
Then
\begin{eqnarray*}
\{ \beta \in R|\ \beta>0 \wedge \mu^{-1}\beta <0 \wedge s_\alpha\beta<0 \}&=&\{ \beta \in R|\ \beta>0 \wedge \mu^{-1}\beta <0 \}\\
&=& \{ e_{2i-1}-e_j |\ 2i<j \} \cup \{ e_{2i}-e_j |\ 2i<j \}
\end{eqnarray*}
hence
\[ \sum_{\substack{\beta \in R^+_\alpha\\\mu^{-1}\beta<0}} (\beta^\vee+s_\alpha\beta^\vee)=-2\sum_{j=2i+1}^n 2e^*_j \in \begin{cases} 4Q^\vee,& 2|n\\ 4Q^\vee+4e^*_n,& 2|n+1 \end{cases}\]
Notice that $2e^*_n \in Q^\vee$, while $e^*_n \notin Q^\vee$. We get
\[ \boxed{\rho(\mu,\alpha) = \alpha^\vee(-1) \cdot 2e^*_n((-1)^n)} \]

%---------------------------------------------------------------------------------
% Section
%=================================================================================

\section{Comparison of splitting invariants} \label{sec:comparison}

Now we would like to employ the results from the previous sections to compare splitting invariants of different tori. More precisely, we'd like to do the
following: Let $(H,s,\eta)$ be an endoscopic triple for $G$, and $T_1,T_2$ be two maximal tori of $G$ that originate from $H$. We'd like to compare the results
of pairing the endoscopic datum $s$ against the splitting invariants for $T_1$ and $T_2$. In order for this to make sense, we need to show that those
invariants reside in a common space. We will show that the endoscopic characters on $H^1(\Gamma,T_1)$ and $H^1(\Gamma,T_2)$ arising from $s$ factor through
certain quotients of these groups and that those quotients can be related.

For a maximal torus $T$ in $G$ put
\begin{eqnarray*}
X_*(T)_{-1}&=&\{ \lambda \in X_*(T)|\ \sigma_T(\lambda) = -\lambda \}\\
IX_*(T)&=&\{ \lambda-\sigma_T(\lambda)|\ \lambda \in X_*(T) \}
\end{eqnarray*}
where $\sigma_T$ is the action of $\sigma$ on $T$. Recall the Tate-Nakayama isomorphism
\[ \frac{X_*(T)_{-1}}{IX_*(T)} = H^{-1}_T(\Gamma,X_*(T)) \rightarrow H^1(\Gamma,T) \]
given by taking cup-product with the canonical class in $H^2(\Gamma,\C^\times)$. Via this isomorphism, the canonical pairing $\hat T \times X_*(T) \rightarrow
\C^\times$ induces a pairing
\[ \hat T \times H^1(\Gamma,T) \rightarrow \C^\times \]
The splitting invariant enters into the construction of the Langlands-Shelstad transfer factors via this pairing.

\begin{lem} \label{lem:emb}
Let $A' \subset A$ be SOS in $R$. Then each element in the canonical set of isomorphisms $T_{A'}^{S_{A\setminus A'}} \rightarrow T_A$ (Definition
\ref{dfn:canonicalset}) induces the same embedding
\[ i_{A',A}: X_*(T_{A'})_{-1} \hookrightarrow X_*(T_A)_{-1} \]
\end{lem}
\pf For an element $\omega \in \Omega$ put
\[ X_*(T_0)_{\omega=-1} = \{ \lambda \in X_*(T_0)|\ \omega(\lambda)=-\lambda \} \]
For any SOS $B$
\[ X_*(T_0)_{S_B=-1} = \tx{span}_\Q(B) \cap X_*(T_0) \]
and any canonical isomorphism $T_0^{S_B} \rightarrow T_B$ identifies $X_*(T_0)_{S_B=-1}$ with $X_*(T_B)_{-1}$ (this identification will of course depend on the
chosen isomorphism).

Fix one canonical isomorphism $T_{A'}^{S_{A\setminus A'}} \rightarrow T_A$. It is the composition of canonical isomorphisms
\[ T_{A'}^{S_{A\setminus A'}} \stackrel{\varphi^{-1}}{\longrightarrow} T_0^{S_A} \stackrel{\psi}{\rightarrow} T_A \]
and hence induces an inclusion as claimed, because $X_*(T_0)_{S_{A'}=-1} \subset X_*(T_0)_{S_A=-1}$. Moreover, any other canonical isomorphism
$T_{A'}^{S_{A\setminus A'}} \rightarrow T_A$ is given by
\[ T_{A'}^{S_{A\setminus A'}} \stackrel{\varphi^{-1}}{\longrightarrow} T_0^{S_A} \stackrel{S_{A''}}{\longrightarrow} T_0^{S_A} \stackrel{\psi}{\rightarrow} T_A \]
for $A'' \subset A \setminus A'$ and clearly $S_{A''}$ acts trivially on $X_*(T_0)_{S_{A'}=-1}$. \qed

The embedding $i_{A',A}$ induces an embedding
\[ \bar i_{A',A}: \frac{X_*(T_{A'})_{-1}}{IX_*(T_{A'})+i_{A',A}^{-1}(IX_*(T_A))} \hookrightarrow \frac{X_*(T_A)_{-1}}{i_{A',A}(IX_*(T_{A'}))+IX_*(T_A)} \]
and via the Tate-Nakayama isomorphism these quotients correspond to quotients of $H^1(\Gamma,T_{A'})$ and $H^1(\Gamma,T_A)$ respectively. We will argue that if
the tori $T_{A'}$ and $T_A$ originate in an endoscopic group $H$, then the endoscopic character factors through these quotients. This provides a means of
comparing the values of the endoscopic character on the cohomology of both tori.

\begin{lem} Let $(H,s,\eta)$ be an endoscopic triple for $G$ and assume that $T_{A'}$ and $T_A$ ($A' \subset A$) originate from $H$, that is, there exist tori $T_1,T_2
\subset H$ and admissible isomorphisms $T_1 \rightarrow T_{A'}$ and $T_2 \rightarrow T_A$. Write $s_{T_{A'}} \in \hat T_{A'}$ and $s_{T_A} \in \hat T_A$ for
the images of $s$ under the duals of these isomorphisms. Assume that there exists a canonical isomorphism $j : T_{A'}^{S_{A\setminus A'}} \rightarrow T_A$ s.t.
$\hat j(s_{T_A}) = s_{T_{A'}}$ (this can always be arranged%
% as long as we are free to choose the admissible embeddings, i.e. modify them by elements of W(R), where W is the corresponding Weyl group
). Then the characters $s_{T_{A'}}$ and $s_{T_A}$ on $H^1(\Gamma,T_{A'})$ resp. $H^1(\Gamma,T_A)$ factor through the above quotients, and the pull-back of the
character $s_{T_A}$ via $\bar i_{A',A}$ equals $s_{T_{A'}}$.
\end{lem}
\pf We identify $H^1(\Gamma,-)$ with $H^{-1}_T(\Gamma,X_*(-))$ via the Tate-Nakayama isomorphism. Because the element $s_{T_A} \in X^*(T_A) \otimes \C^\times$
is $\Gamma$-invariant, its action on $X_*(T_A)$ annihilates the submodule $IX_*(T_A)$. Thus, the action of $j^*(s_{T_A}) \in X^*(T_{A'}) \otimes \C^\times$ on
$X_*(T_{A'})$ annihilates the submodule $j_*^{-1}(IX_*(T_A))$. But we have arranged things so that $j^*(s_{T_A})=s_{T_{A'}}$ and we see that the action of
$s_{T_{A'}}$ on $X_*(T_{A'})$ via the standard pairing annihilates the submodule $IX_*(T_{A'})+j_*^{-1}(IX_*(T_A))$. By the same argument, the action of
$s_{T_A}$ on $X_*(T_A)$ annihilates the submodule $IX_*(T_A) + j_*(IX_*(T_{A'}))$. Finally notice that by Lemma \ref{lem:emb} the restriction of $j_*$ to
$X_*(T_{A'})_{-1}$ coincides with $i_{A',A}$.\qed

%---------------------------------------------------------------------------------
% Bibliography
%=================================================================================

\begin{table}[b]
Tasho Kaletha\\
tkaletha@math.uchicago.edu\\
University of Chicago\\
5734 S. University Avenue\\
Chicago, Illinois 60637
\end{table}

\end{document}